\documentclass[10pt,journal]{IEEEtran}
\usepackage{graphicx}          
\usepackage{amsmath}           
\usepackage{amssymb}
\usepackage{algorithm}
\usepackage{algpseudocode}
\usepackage{hyperref}
\usepackage{color}
\usepackage{cite}
\usepackage{epstopdf}
\usepackage[latin1]{inputenc}  

\usepackage{amsthm}

\theoremstyle{definition}

\theoremstyle{remark}

\theoremstyle{plain}

\DeclareMathOperator{\E}{E}
\DeclareMathOperator*{\argmin}{arg\,min}
\DeclareMathOperator*{\argmax}{arg\,max}
\newcommand{\mbf}[1]{\mathbf{#1}}
\newcommand{\mbs}[1]{\boldsymbol{#1}}
\newcommand{\what}[1]{\widehat{#1}}
\newcommand{\wtilde}[1]{\widetilde{#1}}

\begin{document}

\title{Estimation for the Linear Model with \\ Uncertain Covariance Matrices}
\author{Dave Zachariah, Nafiseh Shariati, Mats Bengtsson, Magnus Jansson and Saikat Chatterjee\thanks{The authors are with the ACCESS Linnaeus Centre, KTH Royal Institute of
Technology, Stockholm. E-mail:
$\{$dave.zachariah, nafiseh, mats.bengtsson, magnus.jansson$\}$@ee.kth.se and saikatchatt@gmail.com. This research has partly been funded by the Swedish Research Council under contracts 621-2011-5847 and 621-2012-4134. The research leading to these results has received funding from the European Research Council under the European Community's Seventh Framework Programme (FP7/2007-2013) / ERC grant agreement n$^\circ$~228044.}}


\maketitle

\begin{abstract}
\textcolor{black}{We derive a maximum a posteriori estimator for the linear
observation model, where the signal and noise covariance matrices
are both uncertain. The uncertainties are treated probabilistically
by modeling the covariance matrices with prior inverse-Wishart
distributions. The nonconvex problem of jointly estimating the signal of interest and the covariance matrices is tackled by a computationally efficient fixed-point iteration as well as an approximate variational Bayes solution. The statistical performance of estimators is compared
numerically to state-of-the-art estimators from the literature and
shown to perform favorably.}
\end{abstract}

\section{Introduction}
The linear observation model
\begin{equation}
\mbf{y} = \mbf{H} \mbf{x} + \mbf{w} \in \mathbb{R}^m,
\label{eq:linearobservation}
\end{equation}
is ubiquitous in \textcolor{black}{signal processing, statistics and
  machine learning,
  cf. \cite{Scharf1991,Kay1993,KailathEtAl2000,Rao1973,Press2005,RaoEtAl2007,Bishop2006}. Applications
  include regression problems, model fitting, functional magnetic
  resonance imaging, finite impulse response
  identification, block data estimation, stochastic channel
  estimation, tracking, sensor fusion and multi-antenna receivers \cite{Kay1993,FristonEtAl1994,Sayed2003,VanTrees&Bell2013}.} Here $\mbf{x} \in \mathbb{R}^n$ denotes the unknown signal of interest, $\mbf{H} \in \mathbb{R}^{m \times n}$ denotes a given matrix with full column rank and $\mbf{w}$ is zero-mean noise. Many estimation procedures rely on prior knowledge of the statistical properties of $\mbf{x}$ and/or $\mbf{w}$. In particular, the covariance matrices $\mbf{P} = \text{Cov}(\mbf{x})$ and $\mbf{R} = \text{Cov}(\mbf{w})$ are assumed to be known. In practice, however, these statistical properties may be subject to uncertainties.
If assigned nominal covariance matrices, $\mbf{P}_0$ and $\mbf{R}_0$, are based on prior knowledge where the statistics are only approximately stationary and/or prior estimates subject to errors, the resulting inaccuracies lead to degradation of estimation performance.

\textcolor{black}{
One approach is to treat the covariance uncertainties deterministically. This entails specifying a class of possible parameter values \cite{Kassam&Poor1985}. For instance, one could model $\mbf{P} = \mbf{P}_0 + \delta \mbf{P}$ and $\mbf{R} = \mbf{R}_0 + \delta \mbf{R}$ and assume that the errors $\delta \mbf{P}$ and $\delta \mbf{R}$ have known bounds on their spectral norms. In this case, \cite{Eldar2006} derived the linear estimator of $\mbf{x}$ that minimizes the worst-case mean square error (MSE) over the specified class of covariance matrices, drawing upon work in \cite{Eldar&Merhav2004,Eldar&Merhav2005}. The problem was shown to be convex and solved in closed form. The `minimax' MSE approach \cite{Verdu&Poor1984}, however, was found to be overly conservative when evaluating its MSE performance. To compensate for this \cite{Eldar2006} also applied a different criterion based on the minimum attainable MSE over the covariance uncertainty class. The `minimax regret' approach aims to minimize the maximum possible deviation from this MSE value. For the problem to be tractable, however, the uncertainty class was restricted such that the eigenvectors of $\mbf{P}$ and $\mbf{R}$ equal the right and left singular vectors of $\mbf{H}$, respectively, and further, that their eigenvalues have known bounds. To circumvent this restriction, \cite{Mittelman&Miller2010} generalized the minimax regret approach and applied it to a wider covariance uncertainty class with element-wise bounds, but only for the signal covariance $\mbf{P}$. Further, unlike \cite{Eldar2006} the resulting estimator is not obtained in closed-form but requires solving a semidefinite program with quartic complexity in signal dimension $n$. In sum, a drawback of the deterministic approaches is the requirement of a restricted parametric class of covariance uncertainties. Further, they are formulated for a single snapshot and do not provide estimates of the signal and noise covariances, both of which are valuable statistical information in certain applications.
}

\textcolor{black}{
A different approach is to treat the covariance uncertainties probabilistically. This entails specifying distributions for the uncertain parameters \cite{Press2005}. For instance, \cite{Tiao&Zellner1964} and
\cite{Svensson&Lundberg2005} model $\mbf{w}$ as a Gaussian random variable and use various prior distributions on $\mbf{R}$. The signal of interest $\mbf{x}$ is modeled with a noninformative prior distribution and therefore no signal covariance matrix $\mbf{P}$ is considered. In \cite{Tiao&Zellner1964}, the prior distribution of $\mbf{R}$ is noninformative resulting in closed-form solutions of the parameter estimates. By contrast, \cite{Svensson&Lundberg2005} consider informative priors for $\mbf{R}$ but require a sampling-based Markov chain Monte Carlo (MCMC) method for solving the problem, which becomes computationally intractable for larger signal dimensions.
}

\textcolor{black}{
In this paper we seek to generalize the probabilistic approach to jointly estimate the signal of interest, as well as the signal \emph{and} noise covariance matrices. Both unknown matrices are modeled as random and independent quantities around the nominal ones, using tractable priors. To the best of the authors' knowledge this has not been addressed and solved in a tractable way in the literature. In this work we use the inverse-Wishart distribution, which is a conjugate prior to the covariance matrix of a Gaussian distribution.  A discussion on the use of this distribution is given in \cite{Press2005,Tiao&Zellner1964, Svensson&Lundberg2005}, where it is shown to be a modified version of the noninformative Jeffreys prior. The inverse-Wishart distribution has also been used in detection problems where the inaccuracies of the nominal covariance matrices arise due to environmental heterogeneity \cite{BidonEtAl2008a, BidonEtAl2008b, WangEtAl2010}.
}

\textcolor{black}{
We show that the maximum a posteriori probability estimator results in a nonconvex optimization problem, but reveals certain connections with the standard estimators. To solve the problem in a computationally efficient manner we formulate a fixed-point iteration. Whilst proving convergence appears intractable, we prove that the iteration does not diverge and illustrate its converge properties empirically. Further, we derive a variational Bayes solution to the problem as a tractable but approximate alternative. Finally, the resulting estimators are evaluated in terms of average performance and robustness.
}

\emph{Notation:} $| \mbf{A} |$ and $\text{tr}\{ \mbf{A} \}$ denote the determinant and trace of $\mbf{A}$, respectively. $\mbf{A} \otimes \mbf{B}$ denotes the Kronecker product of matrices and $\| \cdot \|_F$ denotes the Frobenius norm. $\mbf{E}_{ij}$ is the $ij$th standard basis matrix. \textcolor{black}{$\mathcal{N}(\mbs{\mu}, \mbf{P})$ denotes a Gaussian distribution with mean $\mbs{\mu}$ and covariance $\mbf{P}$. The inverse-Wishart distribution with parameters $\nu$ and $\mbf{C}$ is denoted $\mathcal{W}^{-1}(\mbf{C}, \nu)$.}

\section{Problem formulation}

\textcolor{black}{
For generality we consider a set of $N$ measurements $\{ \mbf{y}_t \}^N_{t=1}$ and corresponding signals of interests $\{ \mbf{x}_t \}^N_{t=1}$. For notational simplicity we write $\mbf{Y} \triangleq \left[ \mbf{y}_1 \cdots \mbf{y}_N \right] \in \mathbb{R}^{m \times N}$ and $\mbf{X} \triangleq \left[ \mbf{x}_1 \cdots \mbf{x}_N \right] \in \mathbb{R}^{n \times N}$. Then the linear observation model \eqref{eq:linearobservation} is written as
\begin{equation}
\mbf{Y} = \mbf{H}\mbf{X} + \mbf{W}.
\end{equation}
It is assumed that the signal and noise follow
independent Gaussian distributions $\mbf{x}_t|\mbf{P} \sim
\mathcal{N}( \mbs{\mu}_t, \mbf{P} )$ and $\mbf{w}_t|\mbf{R} \sim
\mathcal{N}( \mbf{0}, \mbf{R} )$.
}

When the covariance matrices are known, the maximum a posteriori (MAP) estimator of $\mbf{X}$ coincides with the familiar linear minimum MSE estimator,
\begin{equation}
\begin{split}
\widehat{\mbf{X}}_{\text{map}} &= \argmax_{\mbf{X} \in \mathbb{R}^{n \times N}} \; p(\mbf{X} | \mbf{Y}) \\
&= ( \mbf{H}^\top \mbf{R}^{-1}\mbf{H} + \mbf{P}^{-1} )^{-1}(
\mbf{H}^\top \mbf{R}^{-1}\mbf{Y} + \mbf{P}^{-1}\mbf{U} ),
\end{split}
\label{eq:MAP}
\end{equation}
where $\mbf{U} \triangleq \left[ \mbs{\mu}_1 \cdots \mbs{\mu}_N \right] \in \mathbb{R}^{n \times N}$ \cite{KailathEtAl2000}. As the uncertainty or variance of the prior of $\mbf{x}_t$ increases,
by setting $\mbf{P} = \sigma^2_x \mbf{I}_n$ and $\sigma^2_x \rightarrow
\infty$, the estimator coincides with the minimum variance unbiased
(MVU) estimator, $\widehat{\mbf{X}}_{\text{map}} \rightarrow
\widehat{\mbf{X}}_{\text{mvu}} = ( \mbf{H}^\top \mbf{R}^{-1}\mbf{H} )^{-1}
\mbf{H}^\top \mbf{R}^{-1}\mbf{Y}$ \cite{Kay1993}. When $\mbf{P}$ and $\mbf{R}$ are not known precisely they are replaced by nominal matrices, $\mathbf{P}_0$ and $\mathbf{R}_0$.

Henceforth the unknown covariance matrices are modeled as random  and independent quantities around the nominal ones, using inverse-Wishart distributions: $\mathbf{P} \sim \mathcal{W}^{-1}(\mbf{C}_x, \nu_x)$ and
$\mathbf{R} \sim \mathcal{W}^{-1}(\mbf{C}_w, \nu_w)$. Assuming
that $\E [ \mbf{P} ] = \mbf{P}_0$ and $\E [ \mbf{R} ] =
\mbf{R}_0$, we have $\mathbf{C}_x = (\nu_x -n -1) \mathbf{P}_0$
and $\mathbf{C}_w = (\nu_w -m -1) \mathbf{R}_0$.  The degrees
of freedom, $\nu_x > n + 1$ and $\nu_w > m + 1$, control the
certainties of $\mbf{P}$ and $\mbf{R}$. Extensions to the complex Gaussian and inverse-Wishart distributions \cite{Maiwald&Kraus2000} are straight-forward.

The goal is to estimate $\mbf{X}$, $\mbf{P}$ and $\mbf{R}$ from the set of observations $\mbf{Y}$.

\section{The CMAP estimator}
\textcolor{black}{
The maximum a posterior estimator with random covariance matrices, henceforth denoted CMAP, is obtained by solving
\begin{equation}
\min_{\mathbf{X} \in \mathbb{R}^{n
    \times N}, \: \mbf{P} \succ \mbf{0}, \mbf{R} \succ \mbf{0}} \; p(
  \mbf{X},   \mbf{P}, \mbf{R} | \mbf{Y} ),
\label{eq:cmapproblem}
\end{equation}
where $p( \mbf{X}, \mbf{P}, \mbf{R} | \mbf{Y} )$ denotes the joint
posterior probability density function (pdf).} By applying Bayes' rule
and introducing
\begin{equation*}
\begin{split}
J(\mbf{X},\mbf{P},\mbf{R})  &\triangleq \ln p( \mbf{Y} | \mbf{X}, \mbf{P}, \mbf{R} ) + \ln p( \mbf{X}, \mbf{P}, \mbf{R} ) \\
&= \ln p( \mbf{Y} | \mbf{X}, \mbf{R} ) + \ln \left( p( \mbf{X} | \mbf{P} ) p( \mbf{P} ) p( \mbf{R} ) \right) \\
&= J_1(\mbf{X,R}) + J_2(\mbf{X,P}),
\end{split}
\end{equation*}
where $ J_1(\mbf{X,R}) = \left[ \ln p( \mbf{Y} | \mbf{X}, \mbf{R} ) + \ln p(
\mbf{R} ) \right]$  and $ J_2(\mbf{X,P}) =  \left[ \ln  p( \mbf{X} |
\mbf{P} ) + \ln p( \mbf{P} ) \right]$, we can tackle the problem by first solving for $\mbf{R}$ and $\mbf{P}$. Then
\begin{equation}
\widehat{\mbf{X}}_{\text{cmap}} = \argmax_{\mbf{X} \in \mathbb{R}^{n
    \times N}} \left[  \max_{\mbf{R} \succ \mbf{0}, \: \mbf{P} \succ \mbf{0}} J_1(\mbf{X},\mbf{R}) +  J_2(\mbf{X}, \mbf{P}) \right].
\label{eq:CMAP}
\end{equation}
We begin by finding the maximizing $\mbf{R}$ and $\mbf{P}$ below.

\subsection{Concentrated cost function}
\label{sec:concentrated}
Let $\tilde{\mbf{y}}_t \triangleq \mbf{y}_t - \mbf{H}\mbf{x}_t$,
$\widetilde{\mbf{Y}} \triangleq \mbf{Y} - \mbf{H}\mbf{X}$ and $\gamma_w \triangleq \nu_w+m+1+N$, so that
\begin{equation}
\begin{split}
J_1(\mathbf{X,R}) =& \ln p(\widetilde{\mbf{Y}}|\mbf{R}) + \ln p(\mbf{R}) \\
=&\sum^N_{t=1} -\frac{1}{2} \ln |\mbf{R}| - \frac{1}{2} \text{tr}\left \{
  \mbf{R}^{-1} \tilde{\mbf{y}}_t\tilde{\mbf{y}}^\top_t \right \} \\
&- \frac{\nu_w+m+1}{2} \ln |\mbf{R}| - \frac{1}{2} \text{tr}\left \{ \mbf{C}_w \mbf{R}^{-1} \right \}  + K \\
=& -\frac{\gamma_w}{2} \ln |\mbf{R}| - \frac{1}{2} \text{tr}\left \{ (\mbf{C}_w + \widetilde{\mbf{Y}}\widetilde{\mbf{Y}}^\top ) \mbf{R}^{-1} \right \} + K\\
=& \frac{\gamma_w}{2}   \left(  - \ln |\mbf{R}| - \text{tr} \{
  \widetilde{\mbf{R}} \mbf{R}^{-1}  \}  \right)   + K,
\end{split}
\end{equation}
where $K$ denotes an unimportant constant and $\widetilde{\mbf{R}} \triangleq \frac{1}{\gamma_w}(\mbf{C}_w +
\widetilde{\mbf{Y}}\widetilde{\mbf{Y}}^\top ) $. Then
\begin{equation*}
\begin{split}
\widetilde{J}_1(\mbf{X,R}) &\triangleq  - \ln |\mbf{R}| - \text{tr}\left \{ \widetilde{\mbf{R}} \mbf{R}^{-1} \right \} \\
&= - \ln| \widetilde{\mbf{R}} \widetilde{\mbf{R}}^{-1} \mbf{R} | - \text{tr}\left \{  \widetilde{\mbf{R}} \mbf{R}^{-1} \right \} \\
&= - \ln| \widetilde{\mbf{R}} ( \mbf{R}^{-1} \widetilde{\mbf{R}})^{-1} | - \text{tr}\left \{ \widetilde{\mbf{R}} \mbf{R}^{-1} \right \} \\
&= - \ln|\widetilde{\mbf{R}}| + \ln| \mbf{R}^{-1} \widetilde{\mbf{R}} | - \text{tr}\left \{ \mbf{R}^{-1} \widetilde{\mbf{R}} \right \}
\end{split}
\end{equation*}
attains its maximum when $\mbf{R}^{-1}\widetilde{\mbf{R}} =
\mbf{I}_m$, or $\mbf{R}^\star = \frac{1}{\gamma_w}( \mbf{C}_w +
\widetilde{\mbf{Y}}\widetilde{\mbf{Y}}^\top )$. Similarly, let $\widetilde{\mbf{X}} \triangleq \mbf{X} - \mbf{U}$ and $\gamma_x \triangleq \nu_x + n + 1+N$, then $\mbf{P}^\star = \frac{1}{\gamma_x}( \mbf{C}_x + \widetilde{\mbf{X}}\widetilde{\mbf{X}}^\top )$. Note that both $\mbf{P}^\star$ and $\mbf{R}^\star$ are functions of $\mbf{X}$.

Plugging back the solution, and using the matrix determinant lemma, yields
\begin{equation*}
\begin{split}
J_1(\mbf{X},\mbf{R}^\star) &= -\frac{\gamma_w}{2} \ln \left| \frac{1}{\gamma_w}( \mbf{C}_w +
  \widetilde{\mbf{Y}}\widetilde{\mbf{Y}}^\top ) \right| - \frac{\gamma_w}{2} \text{tr} \left \{ \mbf{I}_m \right \} + K \\
&=- \frac{\gamma_w}{2} \ln \left( \frac{1}{\gamma^m_w} \left | \mbf{C}_w + \widetilde{\mbf{Y}}\widetilde{\mbf{Y}}^\top \right| \right)  + K' \\
&= - \frac{\gamma_w}{2} \ln \left| \mbf{C}_w + \widetilde{\mbf{Y}}\widetilde{\mbf{Y}}^\top \right|   + K'' \\
&= - \frac{\gamma_w}{2} \ln \left( |\mbf{C}_w| \left| \mbf{I}_N +
  \widetilde{\mbf{Y}}^\top \mbf{C}^{-1}_w \widetilde{\mbf{Y}} \right| \right) + K'' \\
&= - \frac{\gamma_w}{2} \ln  \left| \mbf{I}_N +
  \widetilde{\mbf{Y}}^\top \mbf{C}^{-1}_w \widetilde{\mbf{Y}} \right| + K'''. \\
\end{split}
\end{equation*}
Similarly,
\begin{equation*}
J_2(\mbf{X},\mbf{P}^\star) = - \frac{\gamma_x}{2} \ln \left | \mbf{I}_N + \widetilde{\mbf{X}}^\top
\mbf{C}^{-1}_x \widetilde{\mbf{X}} \right | + K.
\end{equation*}
In sum, the optimal estimator is given by
\begin{equation}
\begin{split}
\widehat{\mathbf{X}}_{\text{cmap}} &= \argmin_{\mbf{X} \in \mathbb{R}^{n
    \times N}} V(\mathbf{X}) ,
\end{split}
\end{equation}
where the concentrated cost function equals
\begin{equation}
\begin{split}
V(\mbf{X}) \triangleq& \frac{\gamma_w}{2} \ln \left | \mbf{I}_N + (\mbf{Y} - \mbf{H}\mbf{X})^\top
\mbf{C}^{-1}_w (\mbf{Y} - \mbf{H}\mbf{X}) \right | \\
&+ \frac{\gamma_x}{2} \ln \left | \mbf{I}_N + (\mbf{X} - \mbf{U})^\top
\mbf{C}^{-1}_x (\mbf{X} - \mbf{U}) \right |.
\end{split}
\label{eq:costfunction}
\end{equation}

Next, we study the properties of the cost function by writing it as $V(\mbf{X}) = V_1(\mbf{X}) + V_2(\mbf{X})$, where
\begin{equation*}
\begin{split}
V_1(\mbf{X}) &= \frac{\gamma_w}{2} \ln | \mbf{A}(\mbf{X}) |  \\
V_2(\mbf{X}) &= \frac{\gamma_x}{2} \ln | \mbf{B}(\mbf{X}) |,
\end{split}
\end{equation*}
and $\mbf{A}(\mbf{X}) \triangleq \mbf{I}_N +
(\mbf{Y}-\mbf{H}\mbf{X})^\top \mbf{C}^{-1}_w (\mbf{Y}-\mbf{H}\mbf{X})
\succ \mbf{0}$  and $\mbf{B}(\mbf{X}) \triangleq \mbf{I}_N + (\mbf{X}-\mbf{U})^\top \mbf{C}^{-1}_x (\mbf{X}-\mbf{U}) \succ \mbf{0}$. While the inner matrices are quadratic functions of $\mbf{X}$, the log-determinant makes $V_1(\mbf{X})$ and $V_2(\mbf{X})$ nonconvex functions. Their minima, however, provide the key for finding minima of $V(\mbf{X})$.

The minimum of $V_1(\mbf{X})$ is $\what{\mbf{X}}_{\text{mvu}}$ and can be verified by computing the gradient. Using the chain-rule,
\begin{equation*}
\frac{\partial V_1}{\partial x_{it}} = \text{tr}
\left\{  \left( \partial_{A} V_1 \right)^\top  \frac{\partial
    \mbf{A}}{\partial x_{it}} \right\},
\end{equation*}
where the inner derivative equals
\begin{equation*}
\begin{split}
\frac{\partial \mbf{A}}{\partial x_{it}} &= \frac{\partial }{\partial x_{it}} \left ( \mbf{I}_N + (\mbf{Y}-\mbf{HX})^\top \mbf{C}^{-1}_w  (\mbf{Y}-\mbf{HX}) \right ) \\
&= - \mbf{E}^\top_{it} \mbf{H}^\top \mbf{C}^{-1}_w(\mbf{Y}-\mbf{HX}) - (\mbf{Y}-\mbf{HX})^\top \mbf{C}^{-1}_w \mbf{H} \mbf{E}_{it},
\end{split}
\end{equation*}
and the outer derivative is $\partial_{A} V_1 = \frac{\gamma_w}{2} \mbf{A}^{-1}$ due to symmetry. Hence
\begin{equation*}
\begin{split}
\frac{\partial V_1}{\partial x_{it}} &= -\frac{\gamma_w}{2} \text{tr}\left\{ \mbf{A}^{-1} \mbf{E}^\top_{it} \mbf{H}^\top \mbf{C}^{-1}_w(\mbf{Y}-\mbf{HX}) \right\} \\ &\quad -\frac{\gamma_w}{2} \text{tr}\left\{ \mbf{A}^{-1} (\mbf{Y}-\mbf{HX})^\top \mbf{C}^{-1}_w \mbf{H} \mbf{E}_{it} \right\} \\
&= - \gamma_w \text{tr}\left\{ \mbf{H}^\top \mbf{C}^{-1}_w(\mbf{Y}-\mbf{HX}) \mbf{A}^{-1} \mbf{E}_{ti} \right\}
\end{split}
\end{equation*}
and
\begin{equation}
\begin{split}
\partial_X V_1 &= - \gamma_w \mbf{H}^\top \mbf{C}^{-1}_w ( \mbf{Y} - \mbf{H} \mbf{X} ) \mbf{A}^{-1}.
\end{split}
\label{eq:dV1}
\end{equation}
Setting $\partial_X V_1 (\mbf{X}) = \mbf{0}$ and solving for $\mbf{X}$ yields the stationary point $\what{\mbf{X}}_{\text{mvu}}$ since $\mbf{C}_w \propto \mbf{R}_0$. Then $\what{\mbf{X}}_{\text{mvu}}$ is the minimizer of $V_1(\mbf{X})$, since $\ln| \cdot |$ is a monotonically increasing function on the set of positive definite matrices and the quadratic function $\mbf{A}(\mbf{X}) \succeq \mbf{A}(\what{\mbf{X}}_{\text{mvu}})$.

\begin{figure*}
  \begin{center}
    \includegraphics[width=1.70\columnwidth]{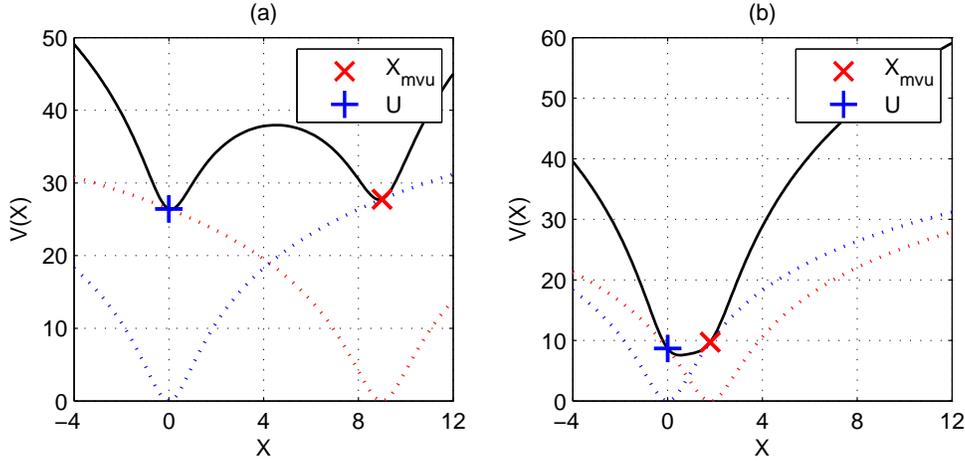}
  \end{center}
  \caption{Example of cost function $V(\mbf{X})$ where $n=1$, $N=1$ and $m=1$ for sake of illustration. Dotted lines show $V_1(\mbf{X})$ and $V_2(\mbf{X})$. Here $\mbf{Y} = \mbf{HX} + \mbf{W}$, where $\mbf{X}=1$ and $\mbf{H}=1$. Nominal variances $\mbf{P}_0 = 0.8$ and $\mbf{R}_0 = 1$ with minimum certainties. (a) $\mbf{W}=8$ resulting in local minima of $V(\mbf{X})$. (b) $\mbf{W}=0.8$ resulting in a single minimum of $V(\mbf{X})$. Note that the minima occur in the vicinity of $\what{\mbf{X}}_{\text{mvu}}$ and $\mbf{U}$.}
  \label{fig:cost_V}
\end{figure*}

Similarly, the trivial minimizer of $V_2(\mbf{X})$ is $\mbf{U}$, and can be verified by
\begin{equation}
\begin{split}
\partial_X V_2 &= \gamma_x \mbf{C}^{-1}_x(\mbf{X}- \mbf{U})\mbf{B}^{-1}.
\end{split}
\label{eq:dV2}
\end{equation}
When a realization $\widehat{\mbf{X}}_{\text{mvu}}$ is far apart from
$\mbf{U}$ then, in the vicinity of the minimizer of $V_1(\mbf{X})$,
$V_2(\mbf{X})$ is approximately constant, and vice versa, due to the
compressive property of the logarithm. In the extreme, therefore,
$V(\mbf{X})$ may have at least two separated minima, located in the vicinity of
$\widehat{\mbf{X}}_{\text{mvu}}$ and $\mbf{U}$, respectively, and the
estimator is not amenable to closed-form solution. On the other hand,
when $\widehat{\mbf{X}}_{\text{mvu}}$ is close to $\mbf{U}$, a single
minimum of $V(\mbf{X})$ may result. These extreme scenarios are illustrated in Fig.~\ref{fig:cost_V}.

Using $\widehat{\mbf{X}}_{\text{mvu}}$ and $\mbf{U}$ as starting
points, minima of $V(\mbf{X})$ can be found by gradient descent
$\what{\mbf{X}}^{\ell+1} = \what{\mbf{X}}^{\ell} - \mu \partial_X V(
\what{\mbf{X}}^{\ell})$, where $\mu > 0$ is the step size and
$\partial_X V = \partial_X V_1 + \partial_X V_2$ given by
\eqref{eq:dV1} and \eqref{eq:dV2}.
The partial derivatives can be written in alternative forms that are
computationally advantageous when $N > n$ and $N > m$, using the matrix inversion lemma,
\begin{equation*}
\begin{split}
 \partial_X V_1 &= - \gamma_w \mbf{H}^\top \mbf{C}^{-1}_w
 \widetilde{\mbf{Y}} ( \mbf{I}_N  +  \widetilde{\mbf{Y}}^\top
 \mbf{C}^{-1}_w \widetilde{\mbf{Y}} )^{-1} \\
&= - \gamma_w \mbf{H}^\top \mbf{C}^{-1}_w
\left( \mbf{I}_N - \widetilde{\mbf{Y}}
  \widetilde{\mbf{Y}}^\top ( \mbf{C}_w +
  \widetilde{\mbf{Y}}\widetilde{\mbf{Y}}^\top ) ^{-1}   \right)
\widetilde{\mbf{Y}} \\
&= - \gamma_w \mbf{H}^\top \mbf{C}^{-1}_w ( \mbf{I}_N +
\widetilde{\mbf{Y}} \widetilde{\mbf{Y}}^\top  \mbf{C}^{-1}_w ) ^{-1}
\widetilde{\mbf{Y}} \\
&= - \gamma_w \mbf{H}^\top ( \mbf{C}_w +
\widetilde{\mbf{Y}} \widetilde{\mbf{Y}}^\top  )^{-1} \widetilde{\mbf{Y}}
\end{split}
\end{equation*}
and similarly
\begin{equation*}
\begin{split}
 \partial_X V_2 &= \gamma_x  \mbf{C}^{-1}_x \widetilde{\mbf{X}}
 (\mbf{I}_N + \widetilde{\mbf{X}}^\top \mbf{C}^{-1}_x
 \widetilde{\mbf{X}} )^{-1} \\
&= \gamma_x ( \mbf{C}_x + \widetilde{\mbf{X}} \widetilde{\mbf{X}}^\top
)^{-1} \widetilde{\mbf{X}}.
\end{split}
\end{equation*}
Thus
\begin{equation*}
\begin{split}
\partial_X V &=  - \gamma_w \mbf{H}^\top \mbf{C}^{-1}_w ( \mbf{Y} - \mbf{H} \mbf{X} ) \mbf{A}^{-1} + \gamma_x \mbf{C}^{-1}_x(\mbf{X}- \mbf{U})\mbf{B}^{-1}\\
&=- \mbf{H}^\top \what{\mbf{R}}^{-1}( \mbf{Y} - \mbf{H} \mbf{X} ) + \what{\mbf{P}}^{-1} (\mbf{X}- \mbf{U}),
\end{split}
\end{equation*}
where
\begin{equation}
\begin{split}
\what{\mbf{P}}(\mbf{X}) &= \frac{1}{\gamma_x} \left( \mbf{C}_x + (\mbf{X}-\mbf{U})(\mbf{X}-\mbf{U})^\top \right) \\
\what{\mbf{R}}(\mbf{X}) &= \frac{1}{\gamma_w} \left( \mbf{C}_w + (\mbf{Y}-\mbf{H}\mbf{X})(\mbf{Y}-\mbf{H}\mbf{X})^\top \right)
\end{split}
\end{equation}
are the covariance matrix estimates. Note that their inverses can be computed recursively by a series of rank-1 updates, using the Sherman-Morrison formula \cite{Hager1989}. The overall computational efficiency of the gradient decent method is, however, dependent on the user-defined step size $\mu$. To circumvent this limitation, we devise an alternative fixed-point iteration method.


\subsection{Fixed-point iteration}

We attempt to find the local minima by iteratively fulfilling the condition for a stationary point. The solution to $\partial_X V(\mbf{X}) = \mbf{0}$, when holding the nonlinear functions $\what{\mbf{P}}(\mbf{X})$ and $\what{\mbf{R}}(\mbf{X})$ constant for a given estimate $\widehat{\mbf{X}}^\ell$, equals
\begin{equation}
\begin{split}
\what{\mbf{X}}^{\ell+1} = ( \mbf{H}^\top \what{\mbf{R}}^{-1}_\ell \mbf{H} + \what{\mbf{P}}^{-1}_\ell )^{-1}(
\mbf{H}^\top \what{\mbf{R}}^{-1}_\ell \mbf{Y} + \what{\mbf{P}}^{-1}_\ell \mbf{U} )
\end{split}
\label{eq:iterate}
\end{equation}
and is iterated until convergence. Comparing \eqref{eq:iterate} with
\eqref{eq:MAP} it is immediately recognized that the fixed-point
method is an iterative application of the standard MAP estimator with
covariance matrices $\what{\mbf{P}}_\ell =
\what{\mbf{P}}(\what{\mbf{X}}^\ell)$ and $\what{\mbf{R}}_\ell = \what{\mbf{R}}(\what{\mbf{X}}^\ell)$. Based on the analysis of the previous section, we propose using $\what{\mbf{X}}_{\text{mvu}}$ and $\mbf{U}$ as two starting points, respectively. The resulting minimum with the lowest cost $V(\mbf{X})$ is then used as the estimate. When the costs happen to be equal, the estimator is indifferent and we can choose the solution that is closest to the MAP estimate, which assumes that the nominal covariances are true. Our numerical experiments show that the iterative solution is very likely to produce the optimal estimate, cf. section~\ref{sec:multipleobservations}.

The CMAP estimator is summarized in Algorithm~\ref{alg:CMAP}. The
function \texttt{iter}$(\cdot)$ iterates \eqref{eq:iterate} until
$\| \widehat{\mbf{X}}^{\ell} - \widehat{\mbf{X}}^{\ell - 1} \|_F < \varepsilon$.

For a derivation of the conditions for convergence of \eqref{eq:iterate} it would be sufficient to prove that the iteration is a contraction mapping \cite{Hasselblatt&Katok2003}. Deriving these conditions appears intractable in general. However, it is possible to show that the iterative solution \eqref{eq:iterate} does not diverge. Let $\what{\mbf{Y}}^{\ell + 1} = \mbf{H}\what{\mbf{X}}^{\ell + 1}$ denote the predicted observation, and $\hat{\mbf{y}}_t$ denote the $t$th column of $\what{\mbf{Y}}^{\ell + 1}$. If $\| \hat{\mbf{y}}_t \|^2_2 = \hat{\mbf{x}}^\top_t \mbf{H}^\top \mbf{H} \hat{\mbf{x}}_t$ is bounded, then $\| \hat{\mbf{x}}_t \|^2_2$ is bounded since $\mbf{H}$ has full rank and $\mbf{H}^\top \mbf{H} \succ \mbf{0}$. Hence $\| \what{\mbf{Y}}^{\ell + 1} \|^2_F < \infty \Rightarrow \| \what{\mbf{X}}^{\ell + 1} \|^2_F < \infty$. Next, consider $\mbf{U}=\mbf{0}$,\footnote{This is no restriction as it is possible to define an equivalent problem with zero-mean variables, $\bar{\mbf{X}} = \mbf{X} - \mbf{U}$ and $\bar{\mbf{Y}} = \mbf{Y} - \mbf{H}\mbf{U}$, and then shift the estimate of $\bar{\mbf{X}}$.} so that \eqref{eq:iterate} can be written as $\what{\mbf{X}}^{\ell+1} = \what{\mbf{P}}_\ell \mbf{H}^\top \left( \what{\mbf{R}}_\ell + \mbf{H} \what{\mbf{P}}_\ell \mbf{H}^\top \right)^{-1}\mbf{Y}$,
and define $\mbs{\Gamma}_\ell \triangleq \mbf{H} \what{\mbf{P}}_\ell \mbf{H}^\top \succ \mbf{0}$ and $\mbs{\Phi}_\ell \triangleq \what{\mbf{R}}_\ell + \mbs{\Gamma}_\ell \succ \mbs{\Gamma}_\ell$. Hence $\| \mbs{\Gamma}_\ell \mbf{\Phi}^{-1}_\ell \|^2_2 < 1$ and it follows that $\| \hat{\mbf{y}}_t \|^2_2 = \| \mbf{H}\hat{\mbf{x}}_t \|^2_2 = \| \mbs{\Gamma}_\ell \mbs{\Phi}^{-1}_\ell \mbf{y}_t \|^2_2 \leq \| \mbs{\Gamma}_\ell \mbs{\Phi}^{-1}_\ell \|^2_2 \| \mbf{y}_t \|^2_2 <  \| \mbf{y}_t \|^2_2$. 
Therefore $\| \what{\mbf{Y}}^{\ell + 1} \|^2_F$ is bounded and consequently $\| \what{\mbf{X}}^{\ell + 1} \|^2_F$ is bounded for all $\ell$. The iterative solution \eqref{eq:iterate} must either converge or produce a bounded orbit. In fact, through extensive simulations the algorithm was always found to converge. In section \ref{sec:empiricalconvergence} we present an empirical convergence analysis of the fixed-point iteration.

\begin{algorithm}
\caption{CMAP estimator} \label{alg:CMAP}
\begin{algorithmic}[1]
\State Input: $\mbf{Y}, \mbf{H}, \mbf{C}_x, \mbf{C}_w, \gamma_x, \gamma_w, \varepsilon$

\State $\widehat{\mbf{X}}^0_1 = \widehat{\mbf{X}}_{\text{mvu}}$ and
$\widehat{\mbf{X}}^0_2 = \mbf{U}$

\State $\widehat{\mbf{X}}_1$ = \texttt{iter}$( \mbf{Y},
\widehat{\mbf{X}}^0_1, \mbf{H}, \mbf{C}_x, \mbf{C}_w, \gamma_x,
\gamma_w, \varepsilon )$

\State $\widehat{\mbf{X}}_2$ = \texttt{iter}$( \mbf{Y},
\widehat{\mbf{X}}^0_2, \mbf{H}, \mbf{C}_x, \mbf{C}_w, \gamma_x,
\gamma_w, \varepsilon )$ \If{$V(\widehat{\mbf{X}}_1) <
V(\widehat{\mbf{X}}_2)$}

\State $\widehat{\mbf{X}} := \widehat{\mbf{X}}_1$
\ElsIf{$V(\widehat{\mbf{X}}_1) > V(\widehat{\mbf{X}}_2)$} \State
$\widehat{\mbf{X}} := \widehat{\mbf{X}}_2$ \Else \State $\widehat{\mbf{X}} :=
\argmin_{ \mbf{X} \in \{ \widehat{\mbf{X}}_1, \widehat{\mbf{X}}_2 \}} \|
\widehat{\mbf{X}}_{\text{map}} - \mbf{X} \|_F$ \EndIf
\State $\widehat{\mbf{P}} = \left(\mbf{C}_x + (\widehat{\mbf{X}}-\mbf{U})(\widehat{\mbf{X}}-\mbf{U})^\top \right)/\gamma_x$
\State $\widehat{\mbf{R}} = \left(\mbf{C}_w + (\mbf{Y}-\mbf{H}\widehat{\mbf{X}})(\mbf{Y}-\mbf{H}\widehat{\mbf{X}})^\top \right)/\gamma_w$
\State Output:
$\widehat{\mbf{X}}$, $\widehat{\mbf{P}}$ and $\widehat{\mbf{R}}$
\end{algorithmic}
\end{algorithm}

\subsection{Marginalized MAP}

In certain applications the covariance matrices $\mbf{P}$ and $\mbf{R}$ may not be of interest and can be treated as nuisance parameters that are marginalized out from the prior and likelihood pdfs. Utilizing the conjugacy of the inverse-Wishart distribution to the Gaussian distribution,
\begin{equation*}
\begin{split}
p(\mbf{X}) &= \int p(\mbf{X} | \mbf{P} ) p(\mbf{P}) d\mbf{P} \\
&\propto \left| \wtilde{\mbf{X}}\wtilde{\mbf{X}}^\top + \mbf{C}_x \right|^{-(\nu_x + N)/2}
\end{split}
\end{equation*}
and
\begin{equation*}
\begin{split}
p(\mbf{Y}|\mbf{X}) &= \int p(\mbf{Y} | \mbf{X}, \mbf{R} ) p(\mbf{R}) d\mbf{R} \\
&\propto \left| \wtilde{\mbf{Y}}\wtilde{\mbf{Y}}^\top + \mbf{C}_w \right|^{-(\nu_w + N)/2}.
\end{split}
\end{equation*}
Then taking the negative logarithm of the marginalized pdf, $p(\mbf{X}|\mbf{Y}) \propto p(\mbf{Y}|\mbf{X}) p(\mbf{X})$, results in a cost function of the same form as \eqref{eq:costfunction} and the marginalized MAP estimator is given by
\begin{equation}
\widehat{\mbf{X}}_{\text{mmap}} = \argmin_{\mathbf{X} \in \mathbb{R}^{n
    \times N}} \; V'(\mbf{X})
\label{eq:MMAP}
\end{equation}
where
\begin{equation*}
\begin{split}
V'(\mbf{X}) \triangleq& \frac{\gamma'_w}{2} \ln \left | \mbf{I}_N + (\mbf{Y} - \mbf{H}\mbf{X})^\top
\mbf{C}^{-1}_w (\mbf{Y} - \mbf{H}\mbf{X}) \right | \\
&+ \frac{\gamma'_x}{2} \ln \left | \mbf{I}_N + (\mbf{X} - \mbf{U})^\top
\mbf{C}^{-1}_x (\mbf{X} - \mbf{U}) \right |,
\end{split}
\end{equation*}
and the weights are $\gamma'_w = \nu_w + N$ and $\gamma'_x = \nu_x + N$. Thus we can apply the same solution methods as used for CMAP but with different weights.

\subsection{Variational MAP}
\textcolor{black}{
We note that the sought variables follow the conditional distributions: $[\mbf{X}]_i| \mbf{P}, \mbf{R}, \mbf{Y} \sim \mathcal{N}( [(\mbf{H}^\top \mbf{R}^{-1} \mbf{H}+ \mbf{P}^{-1})^{-1}(\mbf{H}^\top \mbf{R}^{-1} \mbf{Y} + \mbf{P}^{-1}\mbf{U})]_i, (\mbf{H}^\top \mbf{R}^{-1} \mbf{H}+ \mbf{P}^{-1})^{-1})$, $\mbf{P}| \mbf{X}, \mbf{R}, \mbf{Y} \sim \mathcal{W}^{-1}( \mbf{C}_x + (\mbf{X}-\mbf{U})(\mbf{X}-\mbf{U})^\top, \nu_x + N)$ and $\mbf{R}| \mbf{X}, \mbf{P}, \mbf{Y} \sim \mathcal{W}^{-1}( \mbf{C}_w + (\mbf{Y}-\mbf{H}\mbf{X})(\mbf{Y}-\mbf{H}\mbf{X})^\top, \nu_w + N)$, where $[\mbf{X}]_i$ denotes the $i$th column of $\mbf{X}$.
This enables a numerical computation of the mean of the posterior pdf $p(\mbf{X},\mbf{P},\mbf{R}|\mbf{Y})$ in \eqref{eq:cmapproblem} by means of Markov Chain Monte Carlo methods, e.g., Gibbs sampling \cite{Bishop2006}. Whilst the posterior mean is the MSE-optimal estimate, the dimensionality of the problem requires a very large number of samples for accurate computation, rendering the sampling methods intractable. For completeness we consider a variational approximation of the posterior pdf  \cite{Smidl&Quinn2010}, and derive the corresponding MAP estimator.} The solution to this approximated problem results in an iteration that converges to a local minimum.

The pdf $p(\mbf{X,P,R|Y})$ is approximated by conditionally independent pdfs $q(\mbf{X|Y}) q(\mbf{P|Y}) q(\mbf{R|Y})$. The distributions that minimize the Kullback-Leibler divergence to
$p(\mbf{X,P,R|Y})$ are given by \cite{Bishop2006}
\begin{equation}
\begin{split}
q(\mbf{X|Y}) &\propto e^{\E_{P,R|Y} [ \ln p(\mbf{X,P,R,Y}) ]} \\
q(\mbf{P|Y}) &\propto e^{\E_{X,R|Y} [ \ln p(\mbf{X,P,R,Y}) ]} \\
q(\mbf{R|Y}) &\propto e^{\E_{X,P|Y} [ \ln p(\mbf{X,P,R,Y}) ]}.
\end{split}
\end{equation}
Using the chain rule and introducing $\mbf{V} = \E_{P|Y} [ \mbf{P}^{-1} ]$ and $\mbf{W} = \E_{R|Y} [ \mbf{R}^{-1} ]$ for notational simplicity, we have
\begin{equation*}
\begin{split}
\ln q(\mbf{X|Y}) &= \E_{P,R|Y} [ \ln p(\mbf{X,P,R,Y}) ] + K_1\\
&= \E_{P,R|Y} [ \ln p(\mbf{Y|X,R}) + \ln
p(\mbf{X|P}) ] + K_2 \\
&= -\frac{1}{2} \text{tr}\{(\mbf{Y-HX})^\top \E_{R|Y} [ \mbf{R}^{-1} ] (\mbf{Y-HX}) \} \\
&\quad -\frac{1}{2} \text{tr}\{ (\mbf{ X - U} )^\top  \E_{P|Y} [ \mbf{P}^{-1} ]  ( \mbf{X - U}) \} + K_3 \\
&= -\frac{1}{2} \text{tr}\{\mbf{
Y^\top W Y - Y^\top W H X - X^\top H^\top W Y}  \\
&\quad + \mbf{X^\top H^\top W H X + X^\top V X} \\
&\quad - \mbf{X^\top V U - U^\top V X + U^\top V U} \} + K_4 \\
&= -\frac{1}{2} \text{tr} \{ \mbf{X^\top (H^\top W H + V) X}  \\
&\quad -\mbf{(Y^\top W H + U^\top V)X}  \\
&\quad - \mbf{X^\top (H^\top W Y + VU)} \} + K_5  \\
&=  -\frac{1}{2} \text{tr}\{ (\mbf{X - \wtilde{U}})^\top\wtilde{\mbf{P}}^{-1} (\mbf{X - \wtilde{U}}) \} + K_6,
\end{split}
\end{equation*}
which equals the functional form of $N$ independent Gaussians with mean and
covariance
\begin{equation}
\begin{split}
\wtilde{\mbf{U}} &= (\mbf{H}^\top \E_{R|Y} [ \mbf{R}^{-1} ] \mbf{H}  + \E_{P|Y} [ \mbf{P}^{-1} ])^{-1} \\
&\quad \times (\mbf{H}^\top \E_{R|Y} [ \mbf{R}^{-1} ] \mbf{Y} + \E_{P|Y} [ \mbf{P}^{-1} ]^{-1} \mbf{U}) \\
\wtilde{\mbf{P}} &= (\mbf{H}^\top \E_{R|Y} [ \mbf{R}^{-1} ] \mbf{H}  + \E_{P|Y} [ \mbf{P}^{-1} ])^{-1} .
\end{split}
\end{equation}
The mean and mode coincide and the variational MAP estimator equals
\begin{equation}
\begin{split}
\what{\mbf{X}}_{\text{vmap}} &= \argmax_{\mbf{X} \in \mathbb{R}^{n \times N}} \; q(\mbf{X|Y}) = \wtilde{\mbf{U}}.
\label{eq:vmap}
\end{split}
\end{equation}
Further,
\begin{equation*}
\begin{split}
\ln q(\mbf{P}|\mbf{Y}) &= \E_{X,R|Y}[ \ln p(\mbf{X,P,R,Y}) ] + K_1 \\
&= \E_{X,R|Y}[ \ln p(\mbf{X|P}) + \ln p(\mbf{P}) ] + K_2 \\
&= \E_{X,R|Y}\Bigl[ -\frac{\nu_x + N + n + 1}{2} \ln |\mbf{P}| \\
&\quad- \frac{1}{2} \text{tr}\{
(\mbf{C}_x + \wtilde{\mbf{X}}\wtilde{\mbf{X}}^\top ) \mbf{P}^{-1} \} \Bigr]  + K_3 \\
&= -\frac{\gamma'_x + n + 1}{2} \ln |\mbf{P}| - \frac{1}{2} \text{tr}\{ \wtilde{\mbf{C}}_x \mbf{P}^{-1} \}  + K_4.
\end{split}
\end{equation*}
This is the functional form of an inverse-Wishart with parameters $\gamma'_x = \nu_x + N$ and
\begin{equation}
\begin{split}
\wtilde{\mbf{C}}_x &= \mbf{C}_x + \E_{X,R|Y}[ \mbf{(X-U)(X-U)}^\top ] \\
&= \mbf{C}_x + N\wtilde{\mbf{P}} + N \wtilde{\mbf{U}}\wtilde{\mbf{U}}^\top  - \wtilde{\mbf{U}}\mbf{U}^\top - \mbf{U}\wtilde{\mbf{U}}^\top + \mbf{UU}^\top.\\
\end{split}
\label{eq:tildeC_x}
\end{equation}
Thus $\mbf{P}^{-1}$ follows a Wishart distribution and $\E_{P|Y}[\mbf{P}^{-1}] =
\gamma'_x \wtilde{\mbf{C}}^{-1}_x$. Similarly,
\begin{equation*}
\begin{split}
\ln q(\mbf{R}|\mbf{Y}) &= \E_{X,R|Y}[ \ln p(\mbf{X,P,R,Y}) ] + K_1 \\
&= -\frac{\gamma'_w + m + 1}{2} \ln |\mbf{R}| - \frac{1}{2} \text{tr}\{ \wtilde{\mbf{C}}_w \mbf{R}^{-1} \}  + K_2
\end{split}
\end{equation*}
has the functional form of an inverse-Wishart distribution with parameters $\gamma'_w = \nu_w + N$ and
\begin{equation}
\begin{split}
\wtilde{\mbf{C}}_w &= \mbf{C}_w + \E_{X,P|Y}[ \mbf{(Y-HX)(Y-HX)}^\top ] \\
&= \mbf{C}_w + \mbf{YY}^\top - \mbf{Y}\wtilde{\mbf{U}}^\top \mbf{H}^\top - \mbf{H} \wtilde{\mbf{U}}\mbf{Y}^\top \\
&\quad + N \mbf{H} \wtilde{\mbf{P}} \mbf{H}^\top + N \mbf{H} \wtilde{\mbf{U}}\wtilde{\mbf{U}}^\top \mbf{H}^\top.\\
\end{split}
\label{eq:tildeC_w}
\end{equation}
Thus $\mbf{R}^{-1}$ follows a Wishart distribution and $\E_{R|Y}[\mbf{R}^{-1}] =
\gamma'_w \wtilde{\mbf{C}}^{-1}_w$.

Inserting these results into \eqref{eq:vmap} the variational MAP estimator is computed iteratively as
\begin{equation*}
\begin{split}
\what{\mbf{X}}_{\text{vmap}} &=  (\gamma'_w  \mbf{H}^\top \wtilde{\mbf{C}}^{-1}_w \mbf{H}
+ \gamma'_x \wtilde{\mbf{C}}^{-1}_x)^{-1} \\
&\quad \times  (\gamma'_w \mbf{H}^\top
\wtilde{\mbf{C}}^{-1}_w \mbf{Y} + \gamma'_x \wtilde{\mbf{C}}^{-1}_x \mbf{U}).
\end{split}
\end{equation*}
The parameters $\wtilde{\mbf{C}}_x$ and $\wtilde{\mbf{C}}_w$ are subsequently updated using \eqref{eq:tildeC_x} and \eqref{eq:tildeC_w}. The iteration is initialized by setting the parameters $\wtilde{\mbf{C}}_x = {\mbf{C}}_x$ and $\wtilde{\mbf{C}}_w = {\mbf{C}}_w$. \textcolor{black}{Experimentally we find that using more informative initialization points, i.e., initializing \eqref{eq:tildeC_x} and \eqref{eq:tildeC_w} with $\wtilde{\mbf{U}} = \what{\mbf{X}}_{\text{mvu}}$ and $\mbf{P} = \mbf{P}_0$ produces virtually identical results.}

\section{Experimental results}

In this section we compare the statistical performance of
$\widehat{\mbf{X}}_{\text{cmap}}$ with other estimators using the distribution of normalized squared errors, $\text{NSE} \triangleq \| \mbf{X} - \widehat{\mbf{X}} \|^2_F / \E[ \| \mbf{X} \|^2_F ]$. The expectation is over all random variables. In particular, we will use the normalized mean square error $\text{NMSE} \equiv
\E\left[ \text{NSE} \right]$ and the complementary cumulative
distribution function (ccdf), $\Pr\{ \text{NSE}
> \kappa \}$. The former measures the average performance of the
estimators and the latter quantifies their robustness to noise and
covariance uncertainties.

We also evaluate the NMSE of the covariance matrix estimates $\widehat{\mbf{P}}$ and $\widehat{\mbf{R}}$ in comparison with the nominal matrices $\mbf{P}_0$ and $\mbf{R}_0$.

The statistical measures are estimated by means of Monte Carlo simulations.

\subsection{Estimators}

We compare $\what{\mbf{X}}_{\text{cmap}}$ with $\what{\mbf{X}}_{\text{map}}$ and $\what{\mbf{X}}_{\text{mvu}}$, that use nominal covariance matrices $\mbf{P}_0$ and $\mbf{R}_0$. For CMAP we set the tolerance parameter $\varepsilon = 10^{-6}$.

For comparison of robustness properties with respect to both signal and noise covariance uncertainties, we also apply the state of the art difference regret estimator (DRE) given in \cite{Eldar2006}, $\what{\mbf{X}}_{\text{dre}}$. This estimator assumes that $\mbf{X}$ is zero mean and is derived on assumption that the spatial correlations of the signal and noise are structured by the singular vectors of $\mbf{H}$. Nevertheless, in \cite{Eldar2006} it is suggested that the estimator can be implemented whether or not this correlation structure is satisfied. The covariance uncertainties are treated deterministically as bounds on the eigenvalues of $\mbf{P}_0$, i.e., $l^x_i  \leq  \lambda^x_i  \leq  u^x_i$ for $i=1,\dots,n$, and $\mbf{R}_0$, i.e., $l^w_j  \leq  \lambda^w_j  \leq  u^w_j$ for $j=1,\dots,m$.

The estimator has the form
\begin{equation}
\what{\mbf{X}}_{\text{dre}} = \mbf{D}_x \mbf{H}^{\top} \left(  \mbf{H}\mbf{D}_x\mbf{H}^\top + \mbf{D}_w \right)^{-1} \mbf{Y}.
\end{equation}
The input covariance matrices are set as $\mbf{D}_x = \mbf{V} \mbs{\Delta}_x \mbf{V}^\top$ and $\mbf{D}_w = \mbf{W} \mbs{\Delta}_w \mbf{W}^\top$, where $\mbf{V}$ and $\mbf{W}$ are eigenvector matrices of $\mbf{P}_0$ and $\mbf{R}_0$, respectively. Further, $\mbs{\Delta}_x = \text{diag}(\delta^x_1,\dots, \delta^x_n )$ and $\mbs{\Delta}_w = \text{diag}(\delta^w_1,\dots, \delta^w_m )$, where
\begin{equation*}
\begin{split}
\delta^x_i &= \alpha_i l^x_i + (1-\alpha_i) u^x_i, \quad i=1,\dots,n \\
\delta^w_i &= \alpha_i l^w_i + (1-\alpha_i) u^w_i, \quad i=1,\dots,n
\end{split}
\end{equation*}
and $\delta^w_i = \lambda^w_i$ for all $i=n+1,\dots, m$. Here
\begin{equation}
\alpha_i = \frac{ \sqrt{l^w_i + u^x_i \sigma^2_i} }{ \sqrt{l^w_i + u^x_i \sigma^2_i} + \sqrt{u^w_i + l^x_i \sigma^2_i} },
\end{equation}
where $\sigma_i$ are the singular values of $\mbf{H}$.

Since the covariance uncertainties are treated probabilistically in this work, selecting deterministic bounds on the eigenvalues can only be done heuristically. Here we have selected, $l_i = (1-\nu^0/\nu) \lambda_i$ and $u_i = (1+\nu^0/\nu) \lambda_i$, where $\nu^0$ denotes the minimum integer value of $\nu$, i.e., $\nu^0_x = n+2$ and $\nu^0_w = m+2$. Thus with minimum certainty of the covariances, the lower bound is 0 and upper bound is $2\lambda_i$. As $\nu \rightarrow \infty$, the bounds become tight.

\subsection{Signal setup}
\label{sec:signalsetup}

For sake of illustration, we consider the problem of estimating a stochastic $2
\times 2$ multiple-input multiple output (MIMO) channel $\mbf{A} \in \mathbb{R}^{2 \times 2}$ from observed signals $\mbf{z}_k = \mbf{A} \mbf{s}_k + \mbf{n}_k$. As is common in wireless communications, this is achieved by transmitting a known training sequence $\mbf{S} = [ \mbf{s}_1 \cdots \mbf{s}_K ]$ \cite{Kotecha&Sayeed2004,Bjornson&Ottersten2010}. Collecting $K$ snapshots and vectorizing, the observation is rewritten as $\mbf{y} = \mbf{H} \mbf{x} +
\mbf{w}$, where $\mbf{x} = \text{vec}(\mbf{A})$ and $\mbf{H} = (
\mbf{S}^\top \otimes \mbf{I}_2 )$. The vectorized channel coefficients $\mbf{x}$ and noise $\mbf{w}$ follow  independent, zero-mean, conditionally Gaussian distributions. $\mbf{S}$ is chosen as a deterministic white sequence with constrained power, $\| \mbf{S} \|^2_2 \equiv 10$. We set $\mbf{P}_0 = \frac{1}{n}\mbf{I}_n$, and $\mbf{R}_0 = \sigma^2_w \mbf{I}_m$. The covariance matrices are drawn according to inverse-Wishart distributions. We consider estimating $N$ channel realizations $\mbf{X} \in \mathbb{R}^{n \times N}$ from observations $\mbf{Y} = \mbf{H}\mbf{X} + \mbf{W}  \in \mathbb{R}^{m \times N}$.

The signal to noise ratio,
\begin{equation*}
\text{SNR} \triangleq \frac{\E\left[ \| \mbf{H} \mbf{X} \|^2_F
\right] }{ \E\left[ \| \mbf{W} \|^2_F \right] } = \frac{ \text{tr}
\{ \mbf{H} \mbf{P}_0 \mbf{H}^\top \} }{ \text{tr} \{ \mbf{R}_0  \}
},
\end{equation*}
is varied in the experiments, i.e., setting $\sigma^2_w =
\text{tr}\{ \mbf{H} \mbf{H}^\top \} /( mn \times \text{SNR} )$. We
consider $K = 8$ snapshots so that $m = 16$ and $n=4$. For $m > n$, the resulting low-rank signal structure enables CMAP to estimate parts of both covariances. When $m=n$, the loss of parameter identifiability makes CMAP rely less on the prior signal statistics at higher SNR levels, thus performing closer to the MVU estimator.

Throughout the experiments we ran $10^5$ Monte Carlo simulations for each signal setup.

\subsection{Results for single observation}

\begin{figure}
  \begin{center}
    \includegraphics[width=1.00\columnwidth]{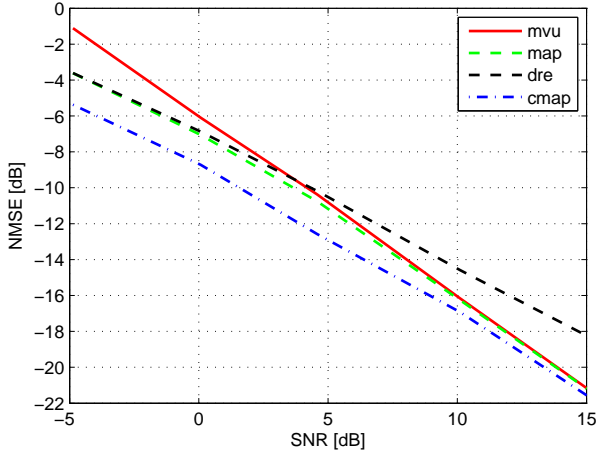}
  \end{center}
  \caption{NMSE versus SNR with minimum certainties of $\mbf{P}$ and $\mbf{R}$. $N=1$.}
  \label{fig:NMSE_N1}
\end{figure}

In the following experiments we consider $N=1$.
First, the average performance of the estimators are compared.
Fig.~\ref{fig:NMSE_N1} shows the NMSE as a function of SNR when the
degrees of freedom for $\mbf{P}$ and $\mbf{R}$ are set to their
minimum integer values, $\nu^0_x = n+2$ and $\nu^0_w = m+2$, respectively.
This yields the minimum certainties of the random quantities.
CMAP is capable of reducing the NMSE by up to approximately 2~dB
compared to the standard MAP. As SNR increases, MAP converges faster
to MVU than does CMAP. The average performance of DRE is initially similar to MAP but the gap increases with SNR as it injects a larger bias.

\begin{figure}
  \begin{center}
    \includegraphics[width=1.00\columnwidth]{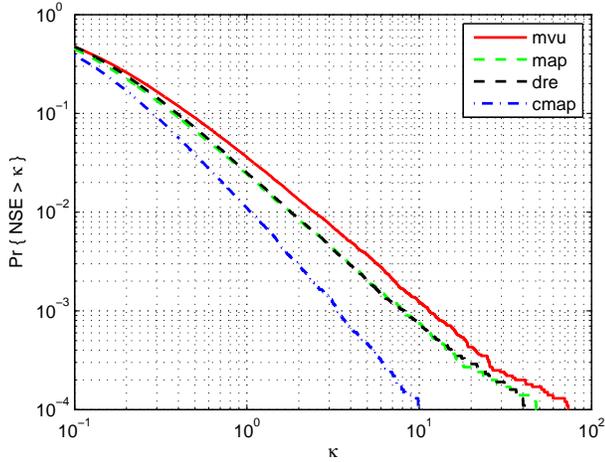}
  \end{center}
  \caption{Ccdf of NSE at SNR=0~dB, with $\nu_x = \nu^0_x$ and $\nu_w = \nu^0_w$. $N=1$.}
  \label{fig:cdf_nux_0_nuw_0_N1}
\end{figure}

Next, the statistical performance of the estimators is compared
using the ccdf, $\Pr \{ \text{NSE} > \kappa \}$, at SNR=0~dB. The curves in Fig.~\ref{fig:cdf_nux_0_nuw_0_N1} illustrate the relative robustness of the estimators to covariance uncertainties. Estimators that produce a lower fraction of poor estimates will have lower ccdfs. Note that $\text{NSE} > 1$ are estimates that have errors greater than the average NSE of the mean, $\what{\mbf{X}} = \mbf{0}$. As expected, MAP and DRE perform similarly at this SNR level, while MVU is slightly worse but declines at a similar rate. CMAP declines more rapidly, with a ccdf that is approximately one order of magnitude lower than MVU at $\kappa = 10$.

\subsection{Results for multiple observations}
\label{sec:multipleobservations}

\begin{figure}
  \begin{center}
    \includegraphics[width=1.00\columnwidth]{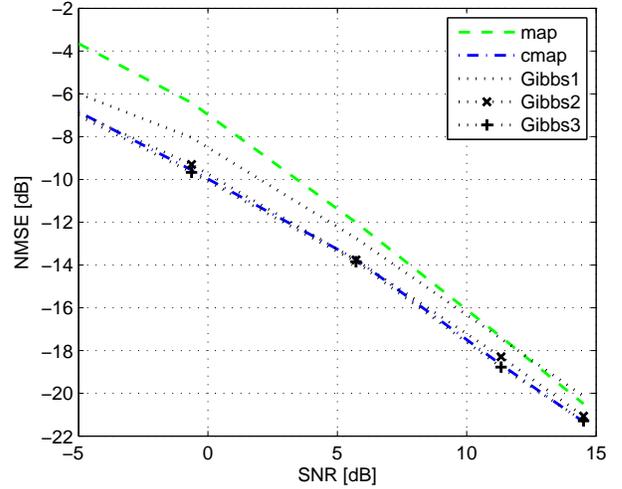}
  \end{center}
  \caption{NMSE versus SNR, with $\nu_x = \nu^0_x$ and $\nu_w = \nu^0_w$. $N=4$. Gibbs 1, 2 and 3 use 200, 2~000 and 20~000 samples, respectively.}
  \label{fig:cdf_nux_0_nuw_0_N4_MCMC}
\end{figure}

\textcolor{black}{The previous experiments are repeated for $N=4$. We
  use the Gibbs sampling approximation of the posterior mean which provides a bound on the NMSE, cf. Fig.~\ref{fig:cdf_nux_0_nuw_0_N4_MCMC}. In this scenario we see that CMAP is very close to the optimum. When $m > N$, the computational complexity of the Gibbs sampler and CMAP is of the order $\mathcal{O}(m^3 N_\text{iter})$, where $m^3$ is the complexity of matrix inversion and $N_\text{iter}$ is the number of repetitions. For the Gibbs sampler, $N_\text{iter}=2\times 10^4$ is about 100 times the number of parameters to estimate and provides a good approximation of the mean. For CMAP, the expected number of iterations is approximately three orders of magnitude less, cf. Sec.~\ref{sec:empiricalconvergence}.}

\begin{figure}
  \begin{center}
    \includegraphics[width=1.00\columnwidth]{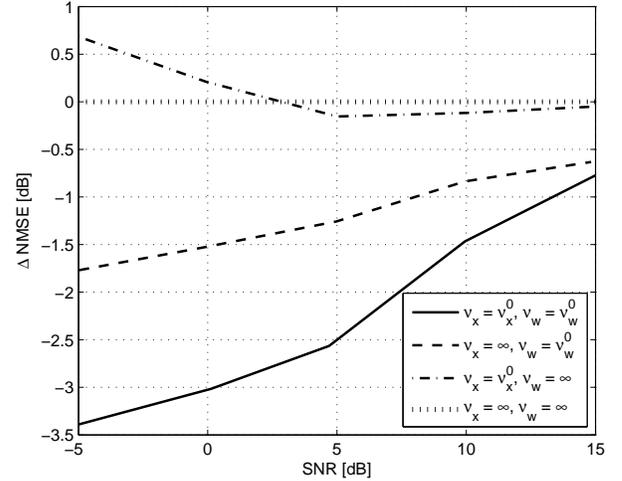}
  \end{center}
  \caption{Difference NMSE between $\what{\mbf{X}}_{\text{cmap}}$ and $\what{\mbf{X}}_{\text{map}}$ versus SNR, for different certainties of $\mbf{P}$ and $\mbf{R}$. $N=4$.}
  \label{fig:delta_NMSE_N4}
\end{figure}

Further, we vary the certainties of the covariance matrices by setting $\nu$ to the extremes, $\nu^0$ and $\infty$. (For $\infty$, we set $\nu$ numerically to $10^5$.) The relative difference in average performance between CMAP and MAP is denoted by $\Delta \text{NMSE}$, where a negative value means reduction in NMSE in decibel using CMAP. The results are shown in Fig.~\ref{fig:delta_NMSE_N4}. When $(\nu_x, \nu_w) = (\infty,\infty)$, CMAP is identical to MAP but for $(\nu_x, \nu_w) = (\nu^0_x, \infty)$ the estimator relies primarily on the noise statistics and CMAP approaches MVU. As both covariance matrices become less certain $(\nu_x, \nu_w) \rightarrow (\nu^0_x, \nu^0_w)$, the advantage of CMAP increases, illustrated by the dashed and solid lines. The improvement for $N=4$ snapshots is above 3~dB for low SNR.

\begin{figure}
  \begin{center}
    \includegraphics[width=1.00\columnwidth]{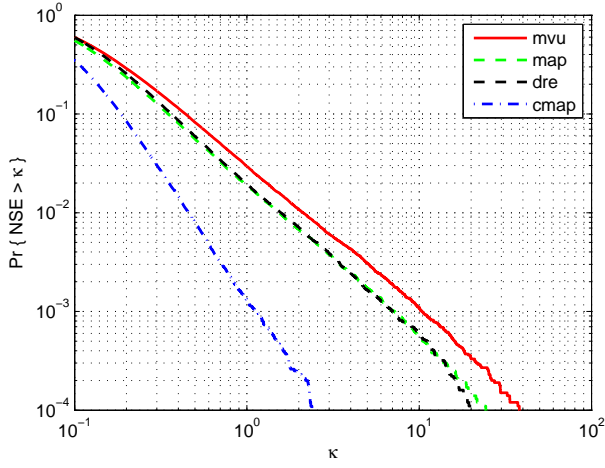}
  \end{center}
  \caption{Ccdf of NSE at SNR=0~dB, with $\nu_x = \nu^0_x$ and $\nu_w = \nu^0_w$. $N=4$.}
  \label{fig:cdf_nux_0_nuw_0_N4}
\end{figure}

\begin{figure}
  \begin{center}
    \includegraphics[width=1.00\columnwidth]{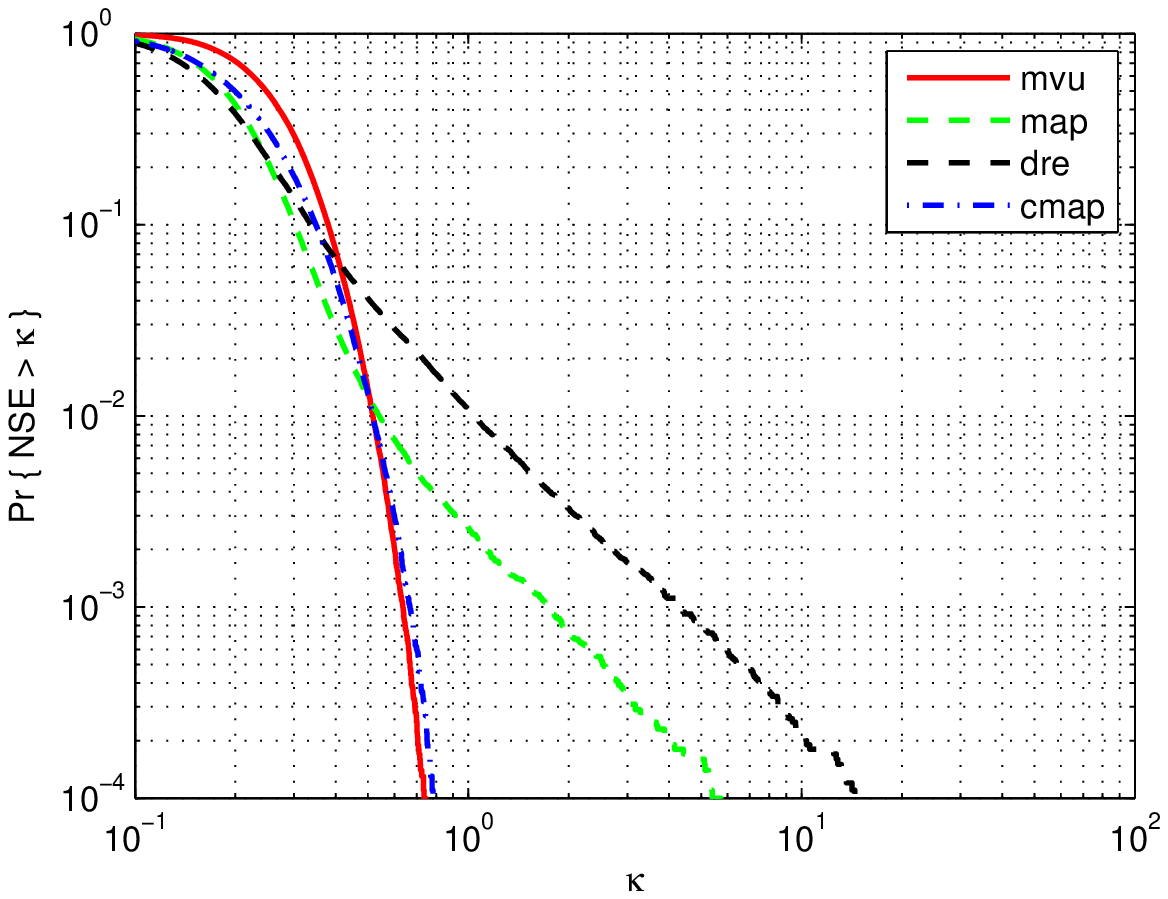}
  \end{center}
  \caption{Ccdf of NSE at SNR=0~dB, with $\nu_x = \nu^0_x$ and $\nu_w = \infty$. $N=4$.}
  \label{fig:cdf_nux_0_nuw_1e5_N4}
\end{figure}

\begin{figure}
  \begin{center}
    \includegraphics[width=1.00\columnwidth]{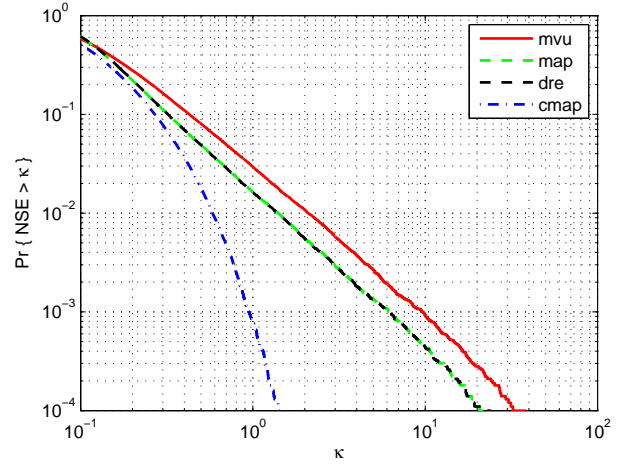}
  \end{center}
  \caption{Ccdf of NSE at SNR=0~dB, with $\nu_x = \infty$ and $\nu_w = \nu^0_w$. $N=4$.}
  \label{fig:cdf_nux_1e5_nuw_0_N4}
\end{figure}

Next, the statistical performance is assessed using the ccdf at SNR=0~dB. Figs.~\ref{fig:cdf_nux_0_nuw_0_N4}, \ref{fig:cdf_nux_0_nuw_1e5_N4} and \ref{fig:cdf_nux_1e5_nuw_0_N4} illustrate robustness at various covariance uncertainties. A comparison between Fig.~\ref{fig:cdf_nux_0_nuw_0_N1} and \ref{fig:cdf_nux_0_nuw_0_N4} shows how the ccdf of CMAP is reduced when the number of samples increases from $N=1$ to $4$. When CMAP relies primarily on the noise statistics, as in Fig.~\ref{fig:cdf_nux_0_nuw_1e5_N4}, it tends towards MVU. While the average NSE of CMAP rises slightly above MAP in this case at low SNR (Fig.~\ref{fig:delta_NMSE_N4}), its ccdf exhibits a sharp decline relative to MAP. When only the signal statistics are reliable, the differences in decline are more pronounced, see Fig.~\ref{fig:cdf_nux_1e5_nuw_0_N4}. In all three cases the fraction of poor estimates cuts off faster for CMAP than MAP.

\begin{figure*}
  \begin{center}
    \includegraphics[width=1.70\columnwidth]{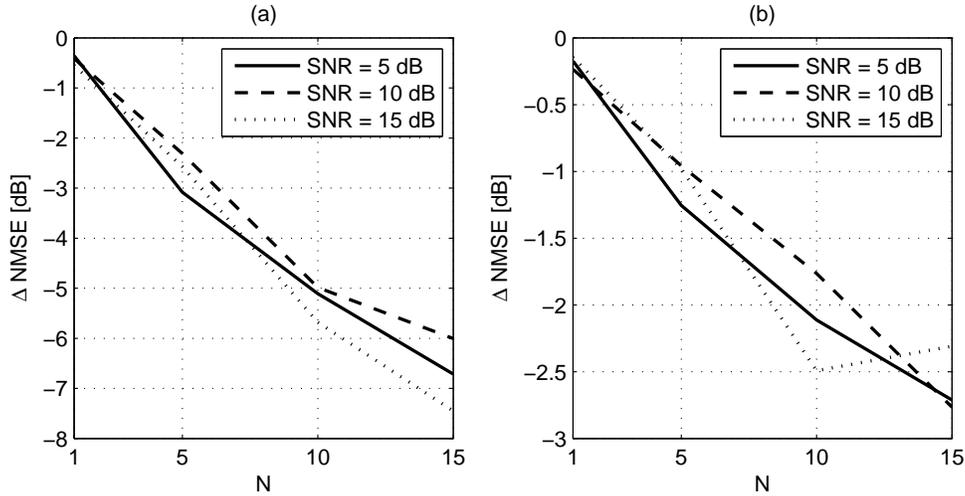}
  \end{center}
  \caption{Difference NMSE between (a) $\what{\mbf{P}}$ and $\mbf{P}_0$ and (b) $\what{\mbf{R}}$ and $\mbf{R}_0$ versus $N$. Here $\nu_x = \nu^0_x$ and $\nu_w = \nu^0_w$.}
  \label{fig:delta_NMSE_cov_N4}
\end{figure*}

Further, we investigate the estimation errors of the covariance matrix estimates $\what{\mbf{P}}$ and $\what{\mbf{R}}$. More specifically, we compute the difference NMSE between the estimates and the priors, which quantifies the information gain, as $N$ increases. The results are shown in Fig.~\ref{fig:delta_NMSE_cov_N4} for various SNR levels. Note that there is a measurable gain even at $N < n$ and $N < m$. Thus CMAP is useful also as a covariance estimator for applications in which signal statistics are of importance.

\begin{figure}
  \begin{center}
    \includegraphics[width=0.90\columnwidth]{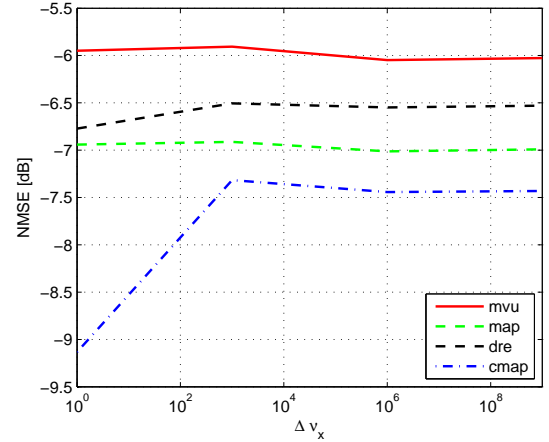}
  \end{center}
  \caption{NMSE versus $\Delta \nu_x$ at SNR=0~dB. $N=4$.}
  \label{fig:nux_N4}
\end{figure}

\begin{figure}
  \begin{center}
    \includegraphics[width=0.90\columnwidth]{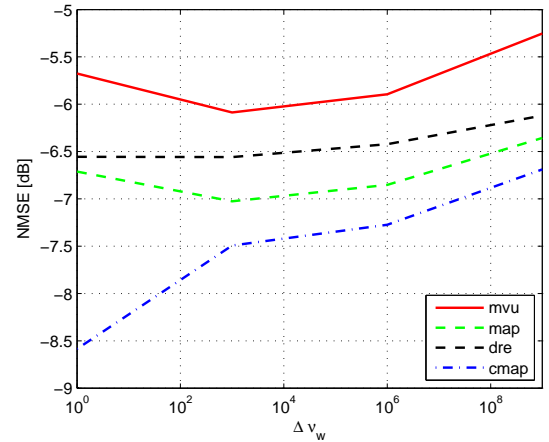}
  \end{center}
  \caption{NMSE versus $\Delta \nu_w$ at SNR=0~dB. $N=4$.}
  \label{fig:nuw_N4}
\end{figure}

In practical scenarios with uncertain covariances,
  CMAP would be implemented with \textcolor{black}{the minimal integer values $\nu^0_x$ and $\nu^0_w$}. We now investigate the effect of mismatches from this conservative prior knowledge by setting the true values to $\nu^0 + \Delta \nu$. In Figures~\ref{fig:nux_N4} and \ref{fig:nuw_N4} we increase $\Delta \nu_x$ and $\Delta \nu_w$, respectively. At SNR=0~dB, we see increases in NMSE for CMAP but the advantage of the estimator is still robust with respect to mismatches for either distributions of $\mbf{P}$ or $\mbf{R}$.

\begin{figure}
  \begin{center}
    \includegraphics[width=1.00\columnwidth]{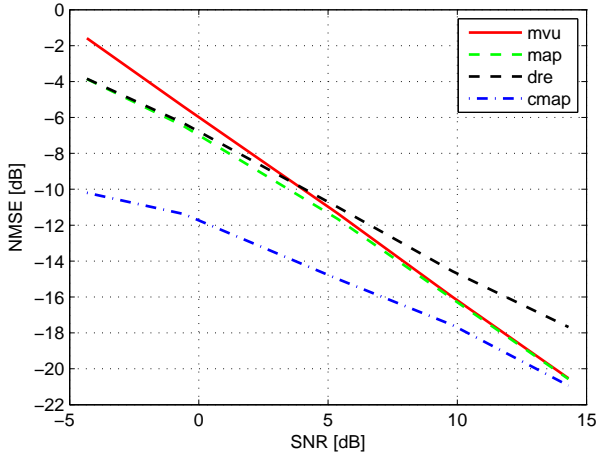}
  \end{center}
  \caption{NMSE versus SNR for $m=64$ and $n=16$. $N=16$.}
  \label{fig:NMSE_N16}
\end{figure}

Finally, we investigate how CMAP performs when increasing the signal dimensions. We now set $m=64$ and $n=16$, for SNR=0~dB, with $\nu_x = \nu^0_x$, $\nu_w = \nu^0_w$ and $N=16$. The NMSE performance as a function of SNR is illustrated in Fig.~\ref{fig:NMSE_N16} which shows a gain of CMAP over MAP greater than in the setup considered in Fig.~\ref{fig:cdf_nux_0_nuw_0_N4_MCMC}, where $m=16$, $m=4$ and $N=4$.

\subsection{Empirical convergence properties}
\label{sec:empiricalconvergence}

\begin{figure}
  \begin{center}
    \includegraphics[width=1.0\columnwidth]{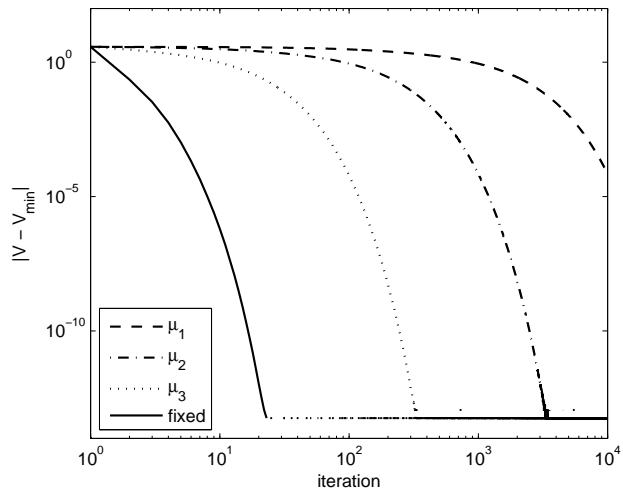}
  \end{center}
  \caption{Convergence to minimum $V_{\text{min}}$ of gradient descent and fixed-point iteration, starting from $\what{\mbf{X}}_{\text{mvu}}$. Step sizes $\mu_1$, $\mu_2$ and $\mu_3$ were set to $10^{-5}$, $10^{-4}$ and $10^{-3}$, respectively. For step-size $10^{-2}$, the gradient descent became unstable. Based on a realization of signal setup given in section \ref{sec:signalsetup} with $N=4$.}
  \label{fig:convergence}
\end{figure}

We now turn to the convergence properties of the iterative solution \eqref{eq:iterate} of CMAP for the same scenario as considered in the previous section, i.e., SNR=0~dB, with $\nu_x = \nu^0_x$, $\nu_w = \nu^0_w$ and $N=4$. Fig.~\ref{fig:convergence} shows a comparison of the convergence rate of the fixed-point iteration and the gradient descent solution, for a typical realization. Both solutions exhibit similar rates once the estimates are sufficiently close to a minimum as both are based on the gradient. But the fixed-point iteration reaches this region within a few iterations without the need for a user-defined step length.

\begin{figure}
  \begin{center}
    \includegraphics[width=1.0\columnwidth]{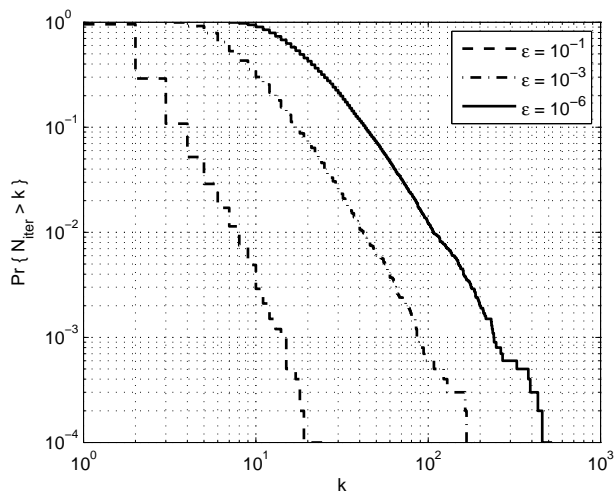}
  \end{center}
  \caption{Ccdf of $N_{\text{iter}}$ at SNR=0~dB, with $\nu_x = \nu^0_x$, $\nu_w = \nu^0_w$ and $N=4$. Tolerances $\varepsilon_1 = 10^{-1}$, $\varepsilon_2 = 10^{-3}$ and $\varepsilon_3 = 10^{-6}$. The corresponding $\E[N_{\text{iter}}]$ was estimated to 2.5,  10.3 and 24.7, respectively, and NMSE was $-9.55$,  $-9.78$ and $-9.78$~dB, respectively.}
  \label{fig:ccdf_iter}
\end{figure}

Next, we study the statistical convergence properties. Let
$N_{\text{iter}}$ denote the total number of iterations until
\eqref{eq:iterate} fulfills $\| \what{\mbf{X}}^{\ell+1} -
\what{\mbf{X}}^{\ell} \|_F < \varepsilon$. Then we can estimate the
ccdf $\Pr \{ N_{\text{iter}} > k \}$, as displayed in
Fig.~\ref{fig:ccdf_iter}. When $\varepsilon = 10^{-6}$ we see that the
probability of $N_{\text{iter}}$ exceeding 300 iterations is less than
$10^{-3}$, and the mean of $N_{\text{iter}}$ is $24.7$. For
$\varepsilon = 10^{-3}$, the mean of $N_{\text{iter}}$ is reduced by
more than a half, while the NMSE is virtually the same. For
$\varepsilon = 10^{-1}$, the entire ccdf is substantially reduced
while incurring an increase in NMSE of only 0.24~dB.

\textcolor{black}{We also estimate the proportion of instances in which the fixed-point iteration converges to two different minima starting from $\what{\mbf{X}}_{\text{mvu}}$ and $\mbf{U}$, respectively. We quantify this as when the convergence points, $\what{\mbf{X}}_1$ and $\what{\mbf{X}}_2$, differ substantially from the numerical tolerance, i.e., $\Pr \{ \| \what{\mbf{X}}_1 - \what{\mbf{X}}_2 \|_F > 10^{-2} \times nN  \}$. For the given scenario, the probability was estimated to 0.03. In 98\% out of those instances $\what{\mbf{X}}_1$ produced a lower cost $V(\mbf{X})$ than $\what{\mbf{X}}_2$.}

\textcolor{black}{In all our simulations we did never encounter a single case when the fixed-point iteration failed to converge to the tolerance.}

\subsection{Comparison between alternative MAP estimators}

\begin{figure}
  \begin{center}
    \includegraphics[width=1.0\columnwidth]{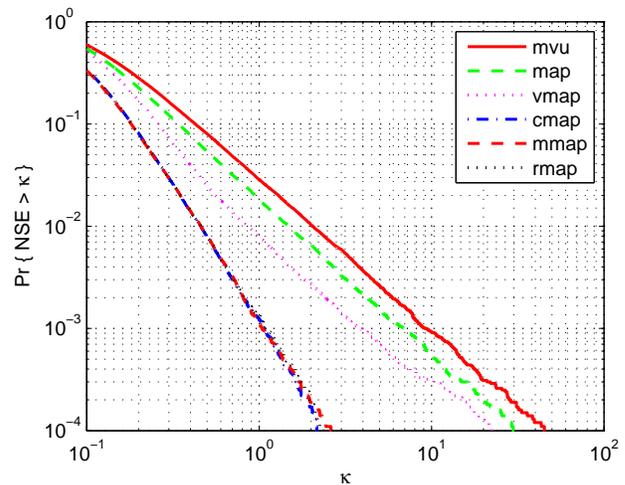}
  \end{center}
  \caption{Ccdf of NSE at SNR=0~dB, with $\nu_x = \nu^0_x$ and $\nu_w = \nu^0_w$. $N=4$. Comparison with the marginalized estimator, `MMAP', the variational estimator, `VMAP', and estimator with randomized initial points `RMAP'.}
  \label{fig:ccdf_mapscomparison}
\end{figure}

Finally, we compare a scenario in which the covariance matrices are not of interest and can be marginalized out, resulting in the marginalized MAP estimator \eqref{eq:MMAP} using the same initial points as CMAP. The difference between the estimators is marginal in terms of NSE performance (see Fig.~\ref{fig:ccdf_mapscomparison}). The NMSE is marginally better for MMAP as it  estimates fewer parameters than CMAP; $-9.98$ and $-9.94$~dB for MMAP and CMAP, respectively. The variational MAP estimator performs better than the standard MAP but is inferior to MMAP and CMAP. The NMSE is $-6.83$ and $-8.10$~dB for MAP and VMAP, respectively.

We also evaluate the significance and robustness of the choice of starting points for CMAP and MMAP. Tests were performed using initial points $\what{\mbf{X}}^0$ randomized by a Gaussian distribution with covariance $\mbf{P}_0$ and a given mean. For each observation $\mbf{Y}$ we then form 10 random initial points $\what{\mbf{X}}^0$ resulting in 10 search paths. The convergence point that yields the lowest cost, $V(\mbf{X})$, is retained as the estimate. We denote this randomized MAP-based estimator as `RMAP'. The different means tested were based on starting points for CMAP, i.e. prior mean $\mbf{U}$ and $\what{\mbf{X}}_{\text{mvu}}$, as well as $\what{\mbf{X}}_{\text{map}}$ and $\what{\mbf{X}}_{\text{dre}}$. Randomizing $\what{\mbf{X}}^0$ around the $\what{\mbf{X}}_{\text{mvu}}$, as well as the proximate values $\what{\mbf{X}}_{\text{map}}$ and $\what{\mbf{X}}_{\text{dre}}$, is found to produce near identical performance to CMAP. Randomizing $\what{\mbf{X}}^0$ around $\mbf{U}$, on the other hand, leads to significantly reduced NSE performance for the worst estimates. These results corroborate the choice of initial points described in section~\ref{sec:concentrated}.
Fig.~\ref{fig:ccdf_mapscomparison} shows the performance for RMAP when using $\what{\mbf{X}}_{\text{dre}}$ as a mean, and the NMSE equals $-9.93$~dB.

\emph{Reproducible research:} Code for reproducing Figs.~\ref{fig:NMSE_N1}, \ref{fig:delta_NMSE_N4} and \ref{fig:NMSE_N16} is available at \url{www.ee.kth.se/~davez/rr-cmap}.

\section{Conclusion}
\textcolor{black}{We have derived a joint signal and covariance maximum a posteriori estimator for the linear observation model, where the signal and noise covariance matrices are modeled as random quantities. We formulated a solution of the nonconvex problem as a fixed-point iterations. The resulting estimator, CMAP, exhibits robustness properties relative to the standard MAP and MVU estimators as well as the minimax difference regret estimator in low-rank signal estimation problems. In this scenario CMAP also shows near MSE-optimal performance. As the number of samples increases, the performance gains of CMAP can be quite substantial.}

\bibliography{refs_cmap}
\bibliographystyle{ieeetr}

\end{document}